# A Robust and Efficient Implementation of LOBPCG

Jed A. Duersch[1], Meiyue Shao[2], Chao Yang[2], and Ming Gu[1,2]

[1]Department of Mathematics, University of California, Berkeley, CA 94720
[2]Computational Research Division, Lawrence Berkeley National Laboratory, Berkeley, CA 94720

April 23, 2018

**Abstract**

Locally Optimal Block Preconditioned Conjugate Gradient (LOBPCG) is widely used to compute eigenvalues of large sparse symmetric matrices. The algorithm can suffer from numerical instability if it is not implemented with care. This is especially problematic when the number of eigenpairs to be computed is relatively large. In this paper we propose an improved basis selection strategy based on earlier work by Hetmaniuk and Lehoucq as well as a robust convergence criterion which is backward stable to enhance the robustness. We also suggest several algorithmic optimizations that improve performance of practical LOBPCG implementations. Numerical examples confirm that our approach consistently and significantly outperforms previous competing approaches in both stability and speed.

# 1 Introduction

Locally Optimal Block Preconditioned Conjugate Gradient (LOBPCG) [9] is a widely used algorithm for computing a few algebraically smallest (or largest) eigenvalues and the corresponding eigenvectors of a Hermitian–definite matrix pencil $(A, B)$ in which $A$ is Hermitian and $B$ is Hermitian positive definite. When $B$ is the identity, the problem becomes a standard symmetric eigenvalue problem.

There are three main advantages of this method when compared with the classical Krylov subspace based methods. First, LOBPCG can employ a good preconditioner when one is available. The use of a preconditioner can dramatically reduce the number of iterations and computation time. Second, the three-term recurrence used by the algorithm keeps memory requirement relatively low, thereby making it possible to tackle problems at a very large scale. Third, because the algorithm is blocked, it can be implemented efficiently on modern parallel computers. Nearly every component of the algorithm can be formulated in level-3 dense or sparse BLAS operations [6, 12] that are highly tuned in several mathematical software libraries (e.g., ACML, Cray LibSci, Intel MKL, etc.) for efficient parallel scaling.

However, it has been observed that the algorithm can breakdown or suffer from numerical instability when it is not implemented carefully. In particular, basis vectors forming the subspace from which the approximate solution to the eigenvalue problem is extracted can become linearly dependent. This problem becomes progressively worse when the number of eigenpairs to be computed becomes relatively large (e.g., hundreds or thousands). For example, in electronic structure calculations, the number of desired eigenpairs is proportional to the number of atoms in the system, which can grow to several thousands [13]. Hence remedies for improving numerical stability are of practical interest.

A strategy proposed in the work of Hetmaniuk and Lehoucq [8] addresses this issue. Their strategy is based on performing additional orthogonalization to ensure that the preconditioned gradient is numerically $B$-orthogonal to both the current and the previous approximations to the desired eigenvectors. However, this strategy can become expensive when the number of eigenpairs to be computed is relatively large. More importantly, reliability can still be severely compromised due to numerical instability within the orthogonalization steps.

This paper presents an efficient and reliable implementation of LOBPCG. We develop a number of techniques to significantly enhance the Hetmaniuk–Lehoucq (HL) orthogonalization strategy in both efficiency and reliability. We also adopt an alternative convergence criterion to ensure achievable error control in computed eigenpairs. For simplicity, we assume that both $A$ and $B$ are real matrices. But our techniques naturally carry over to complex Hermitian matrices.

The rest of this paper is organized as follows. In Section 2, we describe the basic LOBPCG algorithm. In Section 3, we discuss numerical difficulties one may encounter in LOBPCG and the HL strategy for overcoming these difficulties. In Section 4, we present our techniques for improving the HL strategy. In Section 5, we present additional techniques for improving all other aspects of LOBPCG. Finally, in Section 6, we report numerical experimental results to illustrate the effectiveness of our techniques.

## 2 The basic LOBPCG algorithm

We denote the eigenvalues of the symmetric–definite pencil $(A, B)$ arranged in an increasing order by $\lambda_1 \leq \lambda_2 \leq \cdots \leq \lambda_n$. Their corresponding eigenvectors are denoted by $x_1$, $x_2$, …, $x_n$. The first $k \leq n$ eigenvectors and eigenvalues are given by $X = [x_1, x_2, \ldots, x_k]$ and $\Lambda = \operatorname{diag}\{\lambda_1, \lambda_2, \ldots, \lambda_k\}$, respectively, satisfying $AX = BX\Lambda$. It is well known that $X$ is the solution to the trace minimization problem

$$\min_{X^T BX = I} \operatorname{trace}(X^T AX). \tag{1}$$

The LOBPCG algorithm developed by Knyazev seeks to solve (1) by using the updating formula

$$X^{(i+1)} = X^{(i)} C_1^{(i+1)} + X_\perp^{(i)} C_2^{(i+1)},$$

for eigenvector approximation, where $X_\perp^{(i)} = \left[W^{(i)},\ P^{(i)}\right]$. Parenthetical superscript indices indicate the matrix is stored in an array that will be overwritten by subsequent iterations. The block $W^{(i)}$ is the preconditioned gradient of the Lagrangian

$$\mathcal{L}(X, \Lambda) = \frac{1}{2} \operatorname{trace}(X^T AX) - \frac{1}{2} \operatorname{trace}\left[(X^T BX - I)\Lambda\right] \tag{2}$$

associated with (1) at $X^{(i)}$,

$$W^{(i)} = K^{-1}(AX^{(i)} - BX^{(i)}\Theta^{(i)})$$

with $\Theta^{(i)} = X^{(i)T} A X^{(i)}$, where $K$ is any preconditioner. The block $P^{(i)}$ is an aggregated update direction from previous searches recursively defined as

$$P^{(i+1)} = X_{\perp}^{(i)} C_2^{(i+1)},$$

with $P^{(1)}$ being an empty block, i.e., $X_{\perp}^{(1)} = W^{(1)}$.

Coefficient matrices $C_1^{(i+1)}$ and $C_2^{(i+1)}$ are determined at each step of LOBPCG by solving the constrained minimization problem (1) within the subspace $\mathcal{S}^{(i)}$ spanned by $X^{(i)}$, $W^{(i)}$, and $P^{(i)}$. That is,

$$\left(S^{(i)T} A S^{(i)}\right) C^{(i+1)} = \left(S^{(i)T} B S^{(i)}\right) C^{(i+1)} \Theta^{(i+1)}, \tag{3}$$

where $S^{(i)}$ is a matrix whose columns are a basis of $\mathcal{S}^{(i)}$ which is constructed as $S^{(i)} = \left[X^{(i)}, X_{\perp}^{(i)}\right]$ with corresponding

$$C^{(i+1)} = \begin{bmatrix} C_1^{(i+1)} & C_{1\perp}^{(i+1)} \\ C_2^{(i+1)} & C_{2\perp}^{(i+1)} \end{bmatrix}.$$

The leading $k$ columns of $C^{(i+1)}$ form

$$C_x^{(i+1)} = \begin{bmatrix} C_1^{(i+1)} \\ C_2^{(i+1)} \end{bmatrix}$$

which are the components used to compute $X^{(i+1)}$. Remaining columns give the orthogonal complement within the search subspace.

The diagonal matrix $\Theta^{(i+1)}$ in (3) contains approximations to the desired eigenvalues. If $k$ smallest eigenpairs are sought then eigenvalues are sorted in ascending order; otherwise if the largest eigenpairs are sought the order is reversed. Solving the projected eigenvalue problem (3) is often referred to as the *Rayleigh–Ritz procedure*. Combining these steps produces Algorithm 1, the basic LOBPCG algorithm. We leave details of convergence control to Section 4.

# 3 Numerical stability and basis selection

A potential numerical instability issue can occur in Algorithm 1 within the Rayleigh–Ritz procedure shown in Algorithm 2, which is used to solve (3). This is due to the fact that the projection $S^T BS$ can be ill-conditioned or rank deficient, which can happen regardless of the conditioning of $B$.

When $S$ is not $B$-orthonormal, we need to first perform a Cholesky factorization of $S^T BS$ to obtain an upper triangular factor $R$ that is used to transform the generalized eigenvalue problem into a standard eigenvalue problem

$$\left[R^{-T} \left(S^T AS\right) R^{-1}\right] Z = Z\Theta, \tag{4}$$

where $R^T R = S^T BS$. The eigenvectors of pencil $(S^T AS, S^T BS)$ can be recovered by $C = R^{-1} Z$.

When $S^T BS$ is poorly conditioned or numerically singular, Cholesky factorization may fail. Even if the factorization succeeds, $R$ may be poorly conditioned. This poorly conditioned $R$ introduces significant roundoff error in the transformed problem (4) as well as the final transformation $C = R^{-1} Z$. This is often problematic since the near linear dependency of columns in $S$ naturally

**Algorithm 1** The basic LOBPCG algorithm

**Input:**

$X^{(0)}$ is $m \times n_x$ matrix of initial approximate eigenvectors.
$n_v \leq n_x$ is the number of converged eigenvectors requested.
$\tau$ is the threshold used to determine eigenpair convergence.

**Output:**

$X$ is $m \times n_v$ matrix of approximate eigenvectors.
$\Lambda$ is $n_v \times n_v$ diagonal matrix of approximate eigenvalues.

**function** $[X, \Lambda] = \texttt{lobpcgK}(X^{(0)}, n_v, \tau)$
  $[C^{(1)}, \Theta^{(1)}] = \texttt{RayleighRitz}(X^{(0)})$.
  $X^{(1)} = X^{(0)} C^{(1)}$.
  $R^{(1)} = AX^{(1)} - BX^{(1)}\Theta^{(1)}$.
  $P^{(1)} = []$.
  **do** $i = 1, 2, \ldots$
    $W^{(i)} = K^{-1} R^{(i)}$.
    $S^{(i)} = \left[X^{(i)}, W^{(i)}, P^{(i)}\right]$.
    $[C^{(i+1)}, \Theta^{(i+1)}] = \texttt{RayleighRitz}(S^{(i)})$.
    $X^{(i+1)} = S^{(i)} C^{(i+1)}(:, 1:n_x)$.
    $R^{(i+1)} = AX^{(i+1)} - BX^{(i+1)}\Theta^{(i+1)}$.
    $P^{(i+1)} = S^{(i)}(:, n_x+1:\texttt{end}) C^{(i+1)}(n_x+1:\texttt{end}, 1:n_x)$.
    Determine number of converged eigenpairs $n_c$.
  **while** $n_c < n_v$
  Return converged eigenpairs in $X$ and $\Lambda$.
**end function**

**Algorithm 2** Rayleigh–Ritz procedure

**Input:**

$S$ is $m \times n_s$ matrix basis for the search subspace.
*Columns must be linearly independent and well-conditioned with respect to the metric $B$.

**Output:**

$C, \Theta \in \mathbb{R}^{n_s \times n_s}$ that satisfy $C^T(S^T BS)C = I_{n_s}$ and $C^T(S^T AS)C = \Theta$, where $\Theta$ is diagonal.

**function** $[C, \Theta] = \texttt{RayleighRitz}(S)$
  $D = \left(\texttt{diag}(S^T BS)\right)^{-1/2}$.
  Cholesky factorize $R^T R = DS^T BSD$.
  Solve symmetric eigenvalue problem $\left(R^{-T} DS^T ASDR^{-1}\right) Z = Z\Theta$.
  $C = DR^{-1} Z$.
**end function**

emerges when some columns of $X^{(i)}$ become accurate eigenvector approximations. The corresponding columns in both $W^{(i)}$ and $P^{(i)}$ become small in magnitude. They are often sources of potential loss of accuracy and stability. We refer to [5, Section 1.7] for further analysis.

A proper implementation of the LOBPCG algorithm should deflate converged eigenvectors by keeping them in $X^{(i)}$ but exclude corresponding columns from $W^{(i)}$ and $P^{(i)}$. This technique is referred to as *soft locking* [10]. In order to produce reliable results, approximate eigenpairs (Ritz pairs) should be locked in order, i.e., the $(j+1)$st Ritz pairs cannot be locked if the $j$th Ritz pairs does not satisfy the convergence criterion. Some implementations allow out-of-order locking which can lead to slightly better performance. But this approach risks missing desired eigenpairs before termination.

Another technique that helps overcome poor scaling is to normalize each column of $W^{(i)}$ and $P^{(i)}$ before performing the Rayleigh–Ritz procedure. This is equivalent to scaling $S^T BS$ by a diagonal matrix $D$. Note that $R^{-1}$ (or $R^{-T}$) is applied three times in (4) and in $C = R^{-1}Z$. The diagonal scaling step often dramatically reduces the condition number of $R$ and hence improves the numerical stability.

Unfortunately, neither soft locking nor simple diagonal scaling can completely eliminate the numerical instability that potentially leads to a breakdown. It is observed that even with soft locking and diagonal scaling, $S^T BS$ can still become ill-conditioned. When the number of eigenpairs to be computed is relatively large, $S^T BS$ can become ill-conditioned before any approximate eigenvectors in $X^{(i)}$ are sufficiently accurate. This type of failure is quite common and is observed in some of the test cases we present in Section 6. We also provide a concrete example of this phenomenon in Section 4.1.

Hetmaniuk and Lehoucq (HL) proposed in [8] a way to overcome the numerical difficulty associated with ill-conditioning in $S^T BS$. Their basic approach is to keep the $X$, $W$ and $P$ blocks in the subspace $\mathcal{S}$ mutually $B$-orthogonal. They refer to this as a basis selection strategy for $\mathcal{S}$.

Assuming the blocks $X^{(i)}$ and $P^{(i)}$ are $B$-orthonormal already, the basis selection strategy proposed by HL is performed in two steps on each iteration:

1. Before the Rayleigh–Ritz procedure, $W^{(i)}$ is obtained from residuals and then $B$-orthogonalized against both $X^{(i)}$ and $P^{(i)}$. Columns of $W^{(i)}$ are then $B$-orthonormalized.

2. After the Rayleigh–Ritz procedure has been performed, $P^{(i+1)}$ is implicitly $B$-orthogonalized against $X^{(i+1)}$. This is done by forming $C_p^{(i+1)}$ from
$$\begin{bmatrix} 0 \\ C_2^{(i+1)} \end{bmatrix}$$
which is orthogonalized against $C_x^{(i+1)}$ in the metric defined by $S^{(i)T} BS^{(i)}$. The result is then orthonormalized in the same metric producing fully orthonormal blocks $X^{(i+1)} = S^{(i)} C_x^{(i+1)}$ and $P^{(i+1)} = S^{(i)} C_p^{(i+1)}$.

HL use a procedure `ortho`() to carry out both of these orthogonalization steps. For implementation, they refer to the SVQB algorithm developed by Stathopolous and Wu [16]. The procedure `ortho`() operates in two nested loops. In the outer loop a candidate basis is orthogonalized against an existing orthonormal basis, called the external basis, using block classical Gram–Schmidt process. In the inner loop, the remainder is orthonormalized using the singular value decomposition. This is done by a function called `svqb`(). These procedures are outlined in Algorithms 3 and 4, respectively. We altered the original versions of these algorithms slightly to incorporate a metric $M$ which will be used as either $M = B$ or $M = S^T BS$. The original versions only considered the standard Euclidean metric, i.e., $M = I$.

By constructing an orthonormal basis for $S$, HL are able to turn the generalized eigenvalue problem (3) into a standard eigenvalue problem in the Rayleigh–Ritz procedure. Thus Cholesky factorization of the projection $S^T BS$ becomes unnecessary.

**Algorithm 3** Block orthogonalization algorithm proposed by Stathopolous and Wu.

**Input:**

$M$ is $m \times m$ symmetric positive definite metric.

$U^{(\mathtt{in})}$ is $m \times n_u$ candidate basis.

$V$ is $M$-orthonormal external basis.

$\tau_{\mathtt{ortho}} > 0$ is relative orthogonality tolerance.

$\mathtt{itmax}_1$ and $\mathtt{itmax}_2$ are maximum iteration counts for outer and inner loops, respectively.

**Output:**

$U^{(\mathtt{out})}$ is $m \times n_u$ with $M$-orthonormal columns that are $M$-orthogonal to $V$.

$\mathtt{span}([U^{(\mathtt{out})},\ V]) \supseteq \mathtt{span}([U^{(\mathtt{in})},\ V])$.

**function** $U^{(\mathtt{out})} = \mathtt{ortho}(M, U^{(\mathtt{in})}, V, \tau_{\mathtt{ortho}})$
    Set $\tau_{\mathtt{replace}}$.
    **do** $i = 1, 2, \ldots, \mathtt{itmax}_1$
        $U = U - V(V^T M U)$.
        **do** $j = 1, 2, \ldots, \mathtt{itmax}_2$
            $U = \mathtt{svqb}(M, U, \tau_{\mathtt{replace}})$.
        **while** $\frac{\|U^T M U - I_{n_u}\|}{\|MU\|\|U\|} > \tau_{\mathtt{ortho}}$
    **while** $\frac{\|V^T M U\|}{\|MV\|\|U\|} > \tau_{\mathtt{ortho}}$
**end function**

**Algorithm 4** Orthonormalization algorithm using SVD proposed by Stathopolous and Wu.

**Input:**

$M$ is $m \times m$ symmetric positive definite metric.

$U^{(\mathtt{in})}$ is $m \times n_u$.

$\tau_{\mathtt{replace}} > 0$ is tolerance.

**Output:**

$U^{(\mathtt{out})}$ is $m \times n_u$ with $M$-orthonormal columns.

$\mathtt{span}(U^{(\mathtt{out})}) \supseteq \mathtt{span}(U^{(\mathtt{in})})$.

**function** $U^{(\mathtt{out})} = \mathtt{svqb}(M, U^{(\mathtt{in})}, \tau_{\mathtt{replace}})$
    $D = \left(\mathtt{diag}\left(U^T M U\right)\right)^{-1/2}$.
    Solve $\left(D U^T M U D\right) Z = Z\Theta$ for $Z, \Theta$.
    **for all** $\theta_j < \tau_{\mathtt{replace}} \max_i(|\theta_i|)$ **do**
        $\theta_j = \tau_{\mathtt{replace}} \max_i(|\theta_i|)$.
    **end for**
    $U = U D Z \Theta^{-1/2}$.
**end function**

# 4 Stability improvements

## 4.1 Basis truncation

Although the HL basis selection algorithm is plausible and has been demonstrated to work well for a practical problem in [8], its effectiveness hinges on the success of `ortho()` in producing an orthonormal basis $S^{(i)}$.

When source columns in $S^{(i)}$ are nearly linearly dependent, the orthogonalization procedure proposed by Stathopoulos and Wu may not always yield an $S^{(i)}$ that has a sufficiently small condition number. Depending on the implementation, `ortho`() may fail to terminate because the orthogonality error threshold might never be satisfied. This is possible even when $B$ is the identity and becomes more vexing when $B$ is ill-conditioned.

Even if `ortho`() terminates after potentially numerous iterations, the returned basis might be so poorly conditioned that some eigenvalues of $S^TAS$ are spurious. That is, they do not represent an approximation to any eigenvalue of $A$.

The following example illustrates how this problem could occur. Let

$$A = \begin{bmatrix} 3 & 1 & & & & \\ 1 & 3 & 1 & & & \\ & 1 & 3 & 1 & & \\ & & \ddots & \ddots & \ddots & \\ & & & 1 & 3 & 1 \\ & & & & 1 & 3 \end{bmatrix}, \qquad B = K = I, \qquad X^{(0)} = \frac{1}{\sqrt{2}} \begin{bmatrix} 1 & 1 \\ -1 & 1 \\ 0 & 0 \\ \vdots & \vdots \\ 0 & 0 \\ 0 & 0 \end{bmatrix}.$$

Then $\Theta^{(0)} = \operatorname{diag}\{2, 4\}$ and $W^{(0)} = AX^{(0)} - X^{(0)}\Theta^{(0)} = [-e_3,\, e_3]/\sqrt{2}$. The span of $[X^{(0)},\, W^{(0)}]$ is equivalent to $\operatorname{span}([e_1,\, e_2,\, e_3])$. However, a straightforward implementation of the `svqb`() algorithm given in Algorithm 3 may return an output $S$ with four (linearly dependent) columns satisfying $\operatorname{span}(S) = \operatorname{span}([e_1,\, e_2,\, e_3])$. If this basis is used in the subsequent `RayleighRitz`() calculation with the assumption that $S^TBS = I$), then we obtain $\theta_1 = 0$ because $S$ is rank deficient. This is a spurious Ritz value since $A$ is positive definite with smallest eigenvalue larger than 1. The connection between orthogonality error and results of Rayleigh–Ritz is analyzed in more detail in [5, Section 1.7].

There are a few alternative orthogonalization methods, such as the Householder-QR method and Gram–Schmidt process with reorthogonalization [7]. But these alternative methods also have their drawbacks. For instance, the Householder-QR method becomes difficult to apply when $B \neq I$; the Gram–Schmidt process has level 1 or level 2 arithmetic intensity and is not very suitable for high performance computing when the number of vectors to be orthogonalized is relatively large. Since the SVQB algorithm, which is similar to the Cholesky-QR method, has low communication cost [4, 16] and is applicable to both standard and generalized eigenvalue problems, we focus on remedies for SVQB.

To overcome the difficulty in `ortho`(), we truncate basis vectors below roundoff error threshold. The eigenvalue decomposition in `svqb`() gives the form

$$U^TMU = \sum_{i=1}^{n_u} \theta_i z_i z_i^T,$$

where we have sorted eigenvalues $\theta_1 \geq \theta_2 \geq \cdots \geq \theta_{n_u}$. Let $t \leq n_u$ be the number of leading eigenpairs that are above the drop threshold: $\theta_t \geq \theta_1 \cdot \tau_{\texttt{drop}}$, where the drop tolerance $\tau_{\texttt{drop}}$ is a small multiple of the machine precision $\mu_\epsilon$ scaled by the dimension of $U$. The retained basis becomes

$$U^{(\texttt{out})} = U\,[z_1,\, z_2,\, \ldots,\, z_t]\operatorname{diag}\{\theta_1,\, \theta_2,\, \ldots,\, \theta_t\}^{-1/2}.$$

This dropping strategy is mentioned in [16] as an alternative to the standard SVQB algorithm (i.e., Algorithm 4) to handle ill-conditioned inputs, and always produces an output $U^{(\texttt{out})}$ which is

**Algorithm 5** Modified orthogonalization procedure.

**Input:**

$M$ is $m \times m$ symmetric positive definite metric.
$U^{(\texttt{in})}$ is $m \times n_u^{(\texttt{in})}$ candidate basis.
$V$ is $M$-orthonormal external basis.
$\tau_{\texttt{ortho}} > 0$ is relative orthogonality tolerance.

**Output:**

$U^{(\texttt{out})}$ is $m \times n_u^{(\texttt{out})}$, with $n_u^{(\texttt{out})} \leq n_u^{(\texttt{in})}$.
$\texttt{span}([U^{(\texttt{out})},\ V]) \supseteq \texttt{span}([U^{(\texttt{in})},\ V])$.

**function** $U^{(\texttt{out})}$=$\texttt{orthoDrop}(M, U^{(\texttt{in})}, V, \tau_{\texttt{ortho}})$
  Set $\tau_{\texttt{replace}}$ and $\tau_{\texttt{drop}}$.
  **do** $i = 1, 2, 3$
    $U = U - V(V^T MU)$.
    **do** $j = 1, 2, 3$
      **if** $j = 1$ **then**
        $U = \texttt{svqb}(M, U, \tau_{\texttt{replace}})$.
      **else**
        $U = \texttt{svqbDrop}(M, U, \tau_{\texttt{drop}})$.
      **end if**
    **while** $\frac{\|U^T MU - I_{n_u}\|}{\|MU\|\|U\|} > \tau_{\texttt{ortho}}$
  **while** $\frac{\|V^T MU\|}{\|MV\|\|U\|} > \tau_{\texttt{ortho}}$
**end function**

relatively well-conditioned. Then one or two iterations in $\texttt{svqb}()$ suffice for returning an orthonormal basis [16, 19]. Consequently, the robustness and efficiency of the LOBPCG algorithm are improved because the large condition number of the basis is an important source of numerical instability [5, 8].

If maintaining a guaranteed minimum basis dimension is necessary, the basis can be padded with randomly generated $B$-orthogonalized columns. However randomized padding has never been necessary or useful in any of our experiments.

Our modifications are outlined in Algorithms 5 and 6. The implementation we tested still allows one call to $\texttt{svqb}()$ without basis truncation in order to extract potentially useful information. However, subsequent iterations of the inner loop switch to the modified version to ensure a successful exit.

## 4.2 Improved basis selection strategy

The Rayleigh–Ritz procedure always computes orthonormal primitive eigenvectors $Z$ satisfying

$$\left(R^{-T} S^T ASR^{-1}\right) Z = Z\Theta, \qquad Z^T Z = I,$$

where $R$ is the Cholesky factor of $S^T BS$, i.e., $S^T BS = R^T R$. The matrices $Z$ and $\Theta$ can be partitioned conformally as

$$Z = \begin{bmatrix} Z_1 & Z_{1\perp} \\ Z_2 & Z_{2\perp} \end{bmatrix}, \qquad \Theta = \begin{bmatrix} \Theta_x & 0 \\ 0 & \Theta_\perp \end{bmatrix}, \tag{5}$$

**Algorithm 6** SVQB with dropping

**Input:**

$M$ is $m \times m$ symmetric positive definite metric.

$U^{(\texttt{in})}$ is $m \times n_u^{(\texttt{in})}$.

$\tau_{\texttt{drop}} > 0$ is tolerance.

**Output:**

$U^{(\texttt{out})}$ is $m \times n_u^{(\texttt{out})}$, with $n_u^{(\texttt{out})} \le n_u^{(\texttt{in})}$.

$\texttt{span}(U^{(\texttt{out})}) = \texttt{span}(U^{(\texttt{in})})$.

**function** $U^{(\texttt{out})} = \texttt{svqbDrop}(M, U^{(\texttt{in})}, \tau_{\texttt{drop}})$

  $D = \left(\texttt{diag}\left(U^T M U\right)\right)^{-1/2}$.

  Solve $\left(DU^T MUD\right) Z = Z\Theta$ for $Z, \Theta$.

  Determine columns to keep $J = \left\{j : \theta_j > \tau_{\texttt{drop}} \max_i(|\theta_i|)\right\}$.

  $U = UDZ(:, J)\Theta(J, J)^{-1/2}$.

**end function**

so that $Z_1$ and $\Theta_x$ are $k \times k$ matrices. For simplicity, the superscripts are omitted. We also omit diagonal scaling step in `RayleighRitz`(), as $D$ can be absorbed into $R$. To construct the new search directions $[X,\ P] = S[C_x,\ C_p]$, HL's basis selection strategy uses

$$C_x = R^{-1} \begin{bmatrix} Z_1 \\ Z_2 \end{bmatrix}, \qquad C_p = R^{-1} \begin{bmatrix} Q_1 \\ Q_2 \end{bmatrix}, \tag{6}$$

where $[Q_1^T,\ Q_2^T]^T$ is obtained by orthonormalizing

$$\begin{bmatrix} 0 \\ Z_2 \end{bmatrix} - \begin{bmatrix} Z_1 \\ Z_2 \end{bmatrix} \begin{bmatrix} Z_1 \\ Z_2 \end{bmatrix}^T \begin{bmatrix} 0 \\ Z_2 \end{bmatrix}$$

using SVQB.

We propose an improved strategy for constructing $C_p$ based on the following observation. Because $Z$ is an orthogonal matrix, we have

$$\begin{bmatrix} 0 \\ Z_2 \end{bmatrix} - \begin{bmatrix} Z_1 \\ Z_2 \end{bmatrix} \begin{bmatrix} Z_1 \\ Z_2 \end{bmatrix}^T \begin{bmatrix} 0 \\ Z_2 \end{bmatrix} = \begin{bmatrix} Z_{1\perp} \\ Z_{2\perp} \end{bmatrix} \begin{bmatrix} Z_{1\perp} \\ Z_{2\perp} \end{bmatrix}^T \begin{bmatrix} 0 \\ Z_2 \end{bmatrix} = \begin{bmatrix} Z_{1\perp} \\ Z_{2\perp} \end{bmatrix} Z_{2\perp}^T Z_2 = -\begin{bmatrix} Z_{1\perp} \\ Z_{2\perp} \end{bmatrix} Z_{1\perp}^T Z_1.$$

Therefore the blocks $Q_1$ and $Q_2$ required in (6) are completely determined by an orthonormal basis of the columns of $[Z_{1\perp}^T,\ Z_{2\perp}^T]^T Z_{1\perp}^T$. Since $[Z_{1\perp}^T,\ Z_{2\perp}^T]^T$ is already orthonormal, the task of finding $Q_1$ and $Q_2$ simplifies to that of finding an orthonormal basis of the columns of $Z_{1\perp}^T$. This can be achieved by performing an LQ or RQ factorization $Z_{1\perp} = L_{1\perp} Q_{1\perp}$. In practice, a Householder reflection based LQ factorization subroutine (e.g., `xGELQF` in LAPACK [2]) can be applied to guarantee the orthogonality to machine precision even if $Z_{1\perp}^T$ does not have full column rank [7, 18]. Finally, we choose $C_p$ as

$$C_p = R^{-1} \begin{bmatrix} Z_{1\perp} \\ Z_{2\perp} \end{bmatrix} Q_{1\perp}^T.$$

This improved strategy replaces the call to `ortho`() required by the HL strategy with an LQ factorization, which is both cheaper and numerically more stable. Similar to HL's basis selection

**Algorithm 7** Rayleigh–Ritz with improved basis selection

**Input:**

$S$ is an $m \times n_s$ matrix forming a basis for the search subspace.
$n_x$ is the number of extreme eigenpair approximations to return.
$n_c$ is the number of converged eigenpairs from the previous iteration.
$\texttt{useOrtho} = \texttt{true}$ indicates $S^T BS = I$.

**Output:**

$C$ is $n_s \times (2n_x - n_c)$. First $n_x$ columns are $C_x$ followed by $n_x - n_c$ giving $C_p$.
Output satisfies $C^T(S^T BS)C = I$ and $C^T(S^T AS)C = \Theta$.

**function** $[C, \Theta, \texttt{useOrtho}]$=$\texttt{RayleighRitzModified}(S,n_x,n_c,\texttt{useOrtho})$
  **if** $\texttt{useOrtho}$ **then**
    Solve $\left(S^T AS\right) Z = Z\Theta$, where $Z$ is partitioned as in (5).
    LQ factorize $L_{1\perp}Q_{1\perp} = Z_{1\perp}$.

$$C_x = \begin{bmatrix} Z_1 \\ Z_2 \end{bmatrix}, \qquad C_p = \begin{bmatrix} Z_{1\perp} \\ Z_{2\perp} \end{bmatrix} Q_{1\perp}^T.$$

  **else**
    $D = \left(\texttt{diag}(S^T BS)\right)^{-1/2}$
    Cholesky factorize $R^T R = DS^T BSD$.
    **if** $(\texttt{cond}(R) > \tau_{\texttt{skip}})$ set $\texttt{useOrtho} = \texttt{true}$ and exit.
    Solve $\left(R^{-T} DS^T ASDR^{-1}\right) Z = Z\Theta$.
    LQ factorize $L_{1\perp}Q_{1\perp} = Z_{1\perp}$.

$$C_x = DR^{-1} \begin{bmatrix} Z_1 \\ Z_2 \end{bmatrix}, \qquad C_p = DR^{-1} \begin{bmatrix} Z_{1\perp} \\ Z_{2\perp} \end{bmatrix} Q_{1\perp}^T.$$

  **end if**
  Update $\Theta$ to represent partial inner products

$$\Theta = \begin{bmatrix} \Theta_x & 0 \\ 0 & Q_p(:, 1 : n_x - n_c)^T \Theta_\perp Q_p(:, 1 : n_x - n_c) \end{bmatrix}.$$

**end function**

strategy, our improve strategy, outlined in Algorithm 7, requires only $\mathcal{O}(k^3)$ computational work. The benefit is to improve the orthogonality of $[X,\ P]$, and avoid explicit orthogonalization on these vectors, which would require $\mathcal{O}(nk^2)$ work. Our strategy ensures that in each step $[X,\ P]$ always has full column rank (as long as the number of columns here does not exceed that of $A$, which is a plausible assumption in practice). Thus basis truncation is only required when orthonormalizing $W$. A simple rounding error analysis can be found in Section 5.2.

## 4.3 Detecting convergence

In some of our numerical experiments with the HL implementation of LOBPCG, we observed both unexpectedly high and low number of iterations required to reach convergence. One such example, `Andrews`, is a standard eigenvalue problem from the SuiteSparse Matrix Collection (formerly the University of Florida Sparse Matrix Collection).[1] `Andrews` is a 60,000 × 60,000 symmetric matrix with 760,154 nonzero elements. We attempted to find the lowest 400 eigenpairs using a block size of 440 columns. We set the convergence criterion to

$$\frac{\|r_i\|_2}{|\theta_i|\|x_i\|_2} = \frac{\|Ax_i - \theta_i x_i\|_2}{|\theta_i|\|x_i\|_2} \leq \tau, \tag{7}$$

which is the default one provided in the Anasazi package [3]. The tolerance is set to $\tau = 10^{-4}$. This test run was forced to exit without having converged after 1,000 iterations.

Another test problem showed the opposite difficulty with convergence. The matrix pencil `filter2D`, which is also available from the SuiteSparse Matrix Collection, has symmetric sparse matrices $A$ and $B$ of dimension 1,668; $A$ has 10,750 nonzero elements and $B$ is diagonal. In this test case we also sought the lowest 400 eigenvalues using 440 columns per block and the same convergence tolerance, and use the convergence criterion

$$\frac{\|r_i\|_2}{|\theta_i|\|x_i\|_B} = \frac{\|Ax_i - \theta_i Bx_i\|_2}{|\theta_i|\|x_i\|_B} \leq \tau. \tag{8}$$

The algorithm reported convergence immediately after the first iteration. We emphasize that the first iteration is merely `RayleighRitz`() on a random matrix.

These difficulties result from the relative residual computation that is typically used to detect convergence. The criteria (7) and (8) have two problems. The first problem is that (8) lack scaling invariance, i.e., scaling $B$ actually changes the criterion. As a result, convergence detection using such a measure is somewhat arbitrary. In the `filter2D` example, $\|B\| \approx 10^{-10}$ which makes convergence too easy to achieve. Of course, this could be repaired by forcing a uniform scaling of $A$ and $B$ or including $\|A\|$ and $\|B\|$ in the denominator. The second problem with relative residual convergence measures persists even if $B$ is the identity. In fact, for an eigenvalue with small magnitude, it may not be possible to compute a residual in floating point arithmetic that satisfies a criterion of the form (7) or (8). In the `Andrews` example, the smallest eigenvalue is less than $10^{-14}$ in magnitude and $\|A\| \approx 10$. Ill-conditioning of $A$ does not permit the smallest eigenvalue to be known to four digits of precision.

We employ an alternative backward stable convergence criterion:

$$\frac{\|r_i\|_2}{\left(\|A\|_2 + |\theta_i|\|B\|_2\right)\|x_i\|_2} \leq \tau. \tag{9}$$

For large matrices the 2-norms on $A$ and $B$ can be estimated with very little computational cost by using a $k \times m$ Gaussian random matrix $\Omega$ with $k \ll m$. The inequality $\|\Omega A\|_F \leq \|\Omega\|_F\|A\|_2$ (see, e.g., [17, Chapter II, Theorem 3.9]) implies that

$$\|A\|_2^{(\Omega)} \stackrel{\text{def}}{=} \frac{\|\Omega A\|_F}{\|\Omega\|_F} \leq \|A\|_2.$$

This guarantees our convergence criterion is satisfied if

$$\frac{\|r_i\|_2}{\left(\|A\|_2 + |\theta_i|\|B\|_2\right)\|x_i\|_2} \leq \frac{\|r_i\|_2}{\left(\|A\|_2^{(\Omega)} + |\theta_i|\|B\|_2^{(\Omega)}\right)\|x_i\|_2} \leq \tau.$$

[1] https://sparse.tamu.edu

Using this convergence test, `Andrews` converges after performing 56 iterations. Likewise, `filter2D` converges after performing 10 iterations. Finally, we remark that the discussion on convergence criterion is valid not only for LOBPCG, but also for more general Hermitian and non-Hermitian eigensolvers.

# 5 Efficiency improvements

The following efficiency improvements can be safely included in LOBPCG without sacrificing algorithmic stability. Our algorithm will be summarized at the end of this section.

## 5.1 Blockwise matrix multiply and update

In the basis selection stage, the update $X^{(i+1)}$ and the search direction $P^{(i+1)}$ are supposed to overwrite the corresponding parts in $S^{(i)}$. However, their calculation also relies on $S^{(i)}$. Therefore, a practical implementation may look like

$$T \leftarrow S^{(i)} \left[C_x^{(i+1)},\, C_p^{(i+1)}\right], \qquad S^{(i+1)}(:, 1 : n_x + n_p) \leftarrow T. \tag{10}$$

At first glance, this step requires a lot of memory for the workspace $T$.

In our implementation, we partition the matrices $S^{(i)}$, $X^{(i+1)}$, and $P^{(i+1)}$ into chunks of row blocks and update them chunk by chunk, i.e.,

$$S^{(i)} = \begin{bmatrix} S_1^{(i)} \\ S_2^{(i)} \\ \vdots \end{bmatrix}, \qquad T \leftarrow S_j^{(i)} \left[C_x^{(i+1)},\, C_p^{(i+1)}\right], \qquad S_j^{(i+1)}(:, 1 : n_x + n_p) \leftarrow T. \tag{11}$$

By this partitioning the memory requirement for the workspace $T$ is largely reduced. Actually very often there is no need to allocate additional memory for $T$, as the space for storing $AW^{(i)}$, which is not needed anymore in the current iteration, can be reused as workspace here.

In addition, we prefer using such a block row partitioning strategy even if there is enough memory to hold a big workspace. It has been observed that the use of (11) often leads to better performance compared to using (10) by directly calling BLAS especially when the number of desired eigenpairs is moderately small. We refer to [1] for more discussions.

We would like to point out that we store $S^{(i)}$ in a way that improve memory access efficiency. In general, contiguous memory access is preferred from the performance perspective. Although it is natural to arrange $S^{(i)}$ in the order of $\left[X^{(i)},\, W^{(i)},\, P^{(i)}\right]$ so that the memory access pattern of $W^{(i)}$ remains unchanged even in the first two iterations of LOBPCG, we use the order $\left[X^{(i)},\, P^{(i)},\, W^{(i)}\right]$ instead. In addition, the columns of $W^{(i)}$ are stored in the opposite order compared to $X^{(i)}$. With this ordering, deflation always drops the rightmost columns in $S^{(i)}$ so that (11) can be performed with a single BLAS call without additional data movement.

Finally, we remark that the blockwise update strategy is used not only in the basis selection strategy. It is also adopted in the SVQB algorithm to accelerate the orthogonalization.

## 5.2 Implicit product updates

When there is sufficient memory to store $S$, $AS$, and $BS$ if $B \neq I$, HL suggest a possible improvement to efficiency by employing implicit product updates. Given block updates $X^{(i+1)} = S^{(i)}C_x^{(i+1)}$ and $P^{(i+1)} = S^{(i)}C_p^{(i+1)}$, matrix products can be implicitly updated using the same transformations:

$$\begin{aligned} AX^{(i+1)} &= AS^{(i)}C_x^{(i+1)}, \qquad & AP^{(i+1)} &= AS^{(i)}C_p^{(i+1)}, \\ BX^{(i+1)} &= BS^{(i)}C_x^{(i+1)}, \qquad & BP^{(i+1)} &= BS^{(i)}C_p^{(i+1)}, \end{aligned}$$

This is beneficial when direct matrix multiplication (i.e., the application of $A$ or $B$) is expensive, which is often the case in practice. Using a similar technique discussed in Section 5.1, we perform these updates by chunks of row blocks.

We extend this technique to reduce computation of block inner products in the projection matrices

$$S^{(i)T}AS^{(i)} = \begin{bmatrix} X^{(i)T}AX^{(i)} & X^{(i)T}AW^{(i)} & X^{(i)T}AP^{(i)} \\ \cdots & W^{(i)T}AW^{(i)} & W^{(i)T}AP^{(i)} \\ \cdots & \cdots & P^{(i)T}AP^{(i)} \end{bmatrix}$$

and

$$S^{(i)T}BS^{(i)} = \begin{bmatrix} X^{(i)T}BX^{(i)} & X^{(i)T}BW^{(i)} & X^{(i)T}BP^{(i)} \\ \cdots & W^{(i)T}BW^{(i)} & W^{(i)T}BP^{(i)} \\ \cdots & \cdots & P^{(i)T}BP^{(i)} \end{bmatrix}.$$

Most implementations we have seen take advantage of some known structure within each block:

$$X^{(i)T}AX^{(i)} = \Theta^{(i)}, \qquad X^{(i)T}BX^{(i)} = I. \tag{12}$$

In fact, we also have

$$X^{(i)T}AP^{(i)} = 0, \qquad X^{(i)T}BP^{(i)} = 0, \qquad P^{(i)}BP^{(i)} = I, \tag{13}$$

because $P^{(i)}$ is in the orthogonal complement of the previous solution. Implicit updating can also be used to avoid an additional block inner product

$$P^{(i)T}AP^{(i)} = C_p^T \left(S^{(i-1)T}AS^{(i-1)}\right) C_p. \tag{14}$$

Every block involving $W^{(i)}$ must be directly computed unless $B = I$ and the preconditioner is $K = I$. In that case, $W^{(i)}$ is the block of residuals which must be orthogonal to the search subspace from which previous Rayleigh–Ritz solutions were formed including both $X^{(i)}$ and $P^{(i)}$ yielding $X^{(i)T}W^{(i)} = 0$ and $W^{(i)T}P^{(i)} = 0$. Since every block inner product must move $\mathcal{O}(2mn_x)$ data through processors, every direct computation avoided significantly improves performance.

However, implicit updates accumulate roundoff error which can sometimes hinder stability and convergence. Let us consider an implicit update of the form

$$\hat{X} = XT, \qquad A\hat{X} = (AX)T,$$

in which $AX$ is kept in memory and multiplied with $T$ from the right. It can be shown (see [5, Section 1.7]) that the computed result, $\texttt{fl}(\texttt{fl}(AX)T)$, in floating point arithmetic, satisfies

$$\texttt{fl}(\texttt{fl}(AX)T) = AXT + \|A\|\|X\|\|T\|E,$$

where $\|E\|$ is a modest multiple of the machine precision $\mu_\epsilon$ scaled by the dimension of $X$. From the roundoff error analysis, we conclude that implicit updates are in general safe for $B = I$ if $T = C_x$ or $T = C_p$ is used when `useOrtho` is true, because $\|T\|_2 = 1$. But if $B \neq I$, $\|C_x\|_2$ and $\|C_p\|_2$ can be as large as $\|B^{-1/2}\|_2$ with `useOrtho` being true. In this case implicit updates should not be used.

We remark that implicit updates should not be used in `ortho`() and `svqb`() when a nontrivial metric $B$ is involved. When performing

$$U^{(1)} = U^{(0)} - V\left((BV)^T U^{(0)}\right), \qquad U^{(2)} = U^{(1)}\left(DZ\Theta^{-1/2}\right),$$

implicit updates of the form

$$BU^{(1)} = BU^{(0)} - BV\left((BV)^T U^{(0)}\right), \qquad BU^{(2)} = BU^{(1)}\left(DZ\Theta^{-1/2}\right)$$

are valid in exact arithmetic. But this is in general a bad idea, especially for $BU^{(2)}$, because the matrix $DZ\Theta^{-1/2}$ can be very ill-conditioned. We have observed that the use of such updates often leads to failure of termination in the inner loop of `ortho`(), if the candidate basis $U^{(0)}$ has a condition number in the metric $B$ of $10^6$ or higher. Nearly all randomized tests will fail under this circumstance. Indeed, the original purpose for using `ortho`() was to handle such ill-conditioned bases that cause LOBPCG to fail otherwise. Therefore, we do not recommend to use implicit update in `ortho`() and `svqb`().

## 5.3 Skipping orthogonalization

Although constructing an orthonormal basis for $S$ guarantees that `RayleighRitz`() will not fail, the construction process itself can be costly and sometimes unnecessary. The principal extra cost is contained in the call to `ortho`() used to complete block $W$. Even if we assume each step within Algorithm 3 succeeds on the first iteration, the corresponding basis update computations would be

$$W \leftarrow W - [X,\, P]([X,\, P]^T BW), \qquad W \leftarrow W(DZ\Theta^{1/2}).$$

As a result, the memory block containing $W$ must be accessed at least twice.

If `ortho`() is skipped, which forces `RayleighRitz`() to construct and apply Cholesky factors of $S^T BS$, the transformations that would have been performed in `ortho`() are subsumed by updates to $X$ and $P$. As a result, iterations that skip `ortho`() reduce memory movement by more than two full passes over $W$. Furthermore, `svqb`() requires solving an $n_x \times n_x$ eigenvalue problem which is much more time consuming than the corresponding Cholesky decomposition when $n_x$ is not tiny.

In order to take advantage of this possible performance improvement without sacrificing stability, we need to determine when the Cholesky decomposition $R^T R = S^T BS$ becomes unreliable. As $R^{-1}$ or $R^{-T}$ is applied three times, a simple heuristic is to require $\texttt{cond}(R)^{-3} \geq \tau_{\texttt{skip}}$, where $\tau_{\texttt{skip}}$ is a modest multiple of the machine precision.

Knowing this allows us to skip orthogonalization of $W$ initially. When the condition number exceeds the safe threshold we switch to iterations that apply full orthogonalization. Note that our method to construct orthogonal $P^{(i)}$ blocks is not expensive and therefore not worth skipping. After $\texttt{cond}(R)$ passes the safe threshold, subsequent iterations tend to remain above the threshold. As a result, we never attempt to switch back to iterations that skip `ortho`().

Another shortcut to take is to exit `ortho`() early when it is safe to do so. Just as conditioning of $S^T BS$ allows us to identify when `ortho`() may be safely skipped, we can also use the condition number of $W^T BW$ within `svqb`() to predict acceptable orthogonality error. The update

**Algorithm 8** Robust and efficient LOBPCG algorithm

**Input:**
$X^{(0)}$ is $m \times n_x$ initial approximate eigenvectors.
$n_v \leq n_x$ is the number of converged eigenvectors requested.
$\tau$ is the threshold used to determine eigenpair convergence.
**Output:**
$X$ is $m \times n_v$ matrix of approximate eigenvectors.
$\Lambda$ is $n_v \times n_v$ diagonal matrix of approximate eigenvalues.

**function** $[X, \Lambda] = \texttt{lobpcg}(X^{(0)}, n_v, \tau)$
  Set $\tau_{\texttt{ortho}}$.
  $[C^{(1)}, \Theta^{(1)}] = \texttt{RayleighRitz}(X^{(0)})$
  $X^{(1)} = X^{(0)} C^{(1)}$.
  $R^{(1)} = AX^{(1)} - BX^{(1)} \Theta^{(1)}$.
  $n_c = 0$; $\texttt{useOrtho} = \texttt{false}$; $P^{(1)} = []$.
  **do** $i = 1, 2, \ldots$
    **if** $\texttt{useOrtho}$ $W^{(i)} = \texttt{orthoDrop}(B, R^{(i)}, [X^{(i)} \ P^{(i)}], \tau_{\texttt{ortho}})$.
    **else** $W^{(i)} = R^{(i)}$.
    $S^{(i)} = \left[X^{(i)}, P^{(i)}, W^{(i)}\right]$.
    $[C^{(i+1)}, \Theta^{(i+1)}, \texttt{useOrtho}] = \texttt{RayleighRitzModified}(S^{(i)}, n_x, n_c, \texttt{useOrtho})$.
    **if** ($\texttt{useOrtho}$ flips to $\texttt{true}$)
      $W^{(i)} = \texttt{orthoDrop}(B, R^{(i)}, [X^{(i)} \ P^{(i)}], \tau_{\texttt{ortho}})$.
      $[C^{(i+1)}, \Theta^{(i+1)}, \texttt{useOrtho}]$
        $= \texttt{RayleighRitzModified}(S^{(i)}, n_x, n_c, \texttt{useOrtho})$.
    $[X^{(i+1)}, P^{(i+1)}] = S^{(i)}[C_x, C_p]$.
    $R^{(i+1)} = AX^{(i+1)} - BX^{(i+1)} \Theta^{(i+1)}$.
    Determine number of converged eigenpairs $n_c$.
  **while** $n_c < n_v$
  Return converged eigenpairs in $X$ and $\Lambda$.
**end function**

$W^{(2)} = W^{(1)}(DZ\Theta^{-1/2})$ produces relative orthogonality error of magnitude $\mu_\epsilon \texttt{cond}(\Theta^{-1/2})$, therefore we can avoid computing the resulting orthogonality error in $W^{(2)T} B W^{(2)}$ when $\texttt{cond}(\Theta^{-1/2})$ is small. Likewise, relative orthogonality error testing in the outer loop of `ortho()` may be skipped if the `svqb()` transformation was well-conditioned on the *first* iteration of the inner loop. This improvement usually allows us to avoid three block inner products for each call to `ortho()`.

# 6 Numerical examples

In this section we demonstrate the efficiency and robustness of our implementation of Algorithm 8 by several examples performed on the Linux Cluster, Edison, at the National Energy Research Scientific Computing Center (NERSC).[2] Each compute nodes of Edison has two 12-core Intel "Ivy Bridge" processors at 2.4 GHz, and 64 GB DDR3 1866 MHz memory. Our tests are performed on a single compute node of Edison, using 24 cores unless otherwise explicitly stated.

[2] `http://www.nersc.gov/users/computational-systems/edison/`

Algorithm 8 is implemented for shared memory architectures using Fortran and OpenMP, and is compared with Blopex [10] and Anasazi [3]. We also implement Hetmaniuk–Lehoucq's LOBPCG algorithm for comparison. To make the performance comparison more meaningful, we add OpenMP directives to Blopex to parallelize most of loop-level operations including dense linear algebra. As in most cases the time to solution instead of the number of iterations is of practical interest, we mainly report the wall clock time in performance tests.

The test matrices we use are shown in Table 1. Most of these test matrices for standard eigenvalue problems are available from the SuiteSparse Matrix Collection, with the exception of `C60`, which is generated from a MATLAB version of the PARSEC [11] software called RSDFT for a Buckyball molecule. The generalized eigenvalue problems are produced from the SIESTA electronic structure software package [15]. The $B$ matrices in these problems correspond to overlap matrices resulting from using a local atomic basis set to discretize the Kohn–Sham equations. These overlap matrices are mildly ill-conditioned with condition numbers on the order of $10^5$. We only list the numbers of nonzeros in the $A$ matrices for these problems. The numbers of nonzeros in the $B$ matrices are slightly less than that in the corresponding $A$ matrices. In all tests we seek roughly 1% (up to 500) of the algebraically smallest eigenvalues (i.e., the leftmost eigenvalues). This is typically required in Kohn–Sham density functional theory based electronic structure calculations in which the ratio between the number of eigenvalues to be computed and the dimension of the matrix is determined by number of degrees of freedom per atom used to discretize the problem. A highly accurate discretization may require 100 degrees of freedom per atom.

Basis blocks are padded by about 10% of the number of desired eigenpairs. The same convergence criterion (9) is used for all implementations, except for tests performed in Section 6.1. The convergence tolerance is set to $10^{-4}$. We have altered the source code of both Blopex and Anasazi implementations of LOBPCG to use our convergence criterion. Programs are allowed to perform up to 2,000 iterations, and up to four wall clock hours.

Table 1: A list of test matrices.

| case | matrix | type | size | nnz | eigenpairs | block dim. |
|---|---|---|---|---|---|---|
| 1 | `C60` | standard | 17,576 | 407,204 | 176 | 194 |
| 2 | `Si5H12` | standard | 19,896 | 738,598 | 199 | 219 |
| 3 | `c-65` | standard | 48,066 | 360,428 | 481 | 529 |
| 4 | `Andrews` | standard | 60,000 | 760,154 | 500 | 550 |
| 5 | `Ga3As3H12` | standard | 61,349 | 5,970,947 | 500 | 550 |
| 6 | `Ga10As10H30` | standard | 113,081 | 6,115,633 | 500 | 550 |
| 7 | nanotube | generalized | 9,984 | 5,076,862 | 105 | 110 |
| 8 | graphene | generalized | 21,060 | 2,357,415 | 203 | 220 |

## 6.1 Convergence comparisons

The number of iterations performed by an implementation of LOBPCG will depend on the method used to detect convergence and the corresponding threshold. This introduces a difficulty when we attempt to compare different implementations. Because we use the modified convergence criterion (9), it is possible that our algorithm benefits from performing fewer iterations compared to other implementations, which obfuscates the meaning of direct timing comparisons.

In order to show that our convergence criterion is a fair improvement, we first compare each implementation using its original convergence method (marked with suffix -OC) and our modification. We also test our own conservative implementation of the Hetmaniuk–Lehoucq algorithm (CHL). This version does not perform any implicit updates (i.e., the techniques discussed in Section 5.2) or employ other performance improvements marked as optional in [8].

The test results are collected in Table 2. For `C60`, Blopex-OC fails due to a negative pivot within Cholesky factorization. In the case of `Andrews`, neither Anasazi-OC nor CHL-OC converges after running for four hours, and hence these results are considered as failure. Since Anasazi-OC and CHL-OC force convergence in order, they cannot terminate as the minimum eigenvalue is too close to zero to satisfy the relative residual threshold in floating point arithmetic. In fact Blopex-OC also does not fully succeed with this example even though it returns converged eigenpairs within the allowed number of iterations, because some of the desired eigenpairs are missing. The same problem arises in the test for `graphene`. Both Blopex-OC and Blopex miss some of the desired eigenpairs. We marked these cases by $*$ in Table 2.

With the exception of Anasazi for `graphene`, we can see from Table 2 that almost all implementations of the LOBPCG algorithm benefit from the new convergence criterion.

Table 2: Execution time for different convergence criteria. Failures (including incorrect solutions marked with $*$) are highlighted in boldface.

| matrix | Anasazi-OC | Anasazi | Blopex-OC | Blopex | CHL-OC | CHL |
|---|---|---|---|---|---|---|
| `C60` | 947.4 | 21.0 | **failed** | 14.2 | 26.6 | 17.2 |
| `Si5H12` | 53.5 | 42.5 | 49.8 | 29.6 | 42.7 | 31.7 |
| `Andrews` | **failed** | 614.1 | **383.7*** | 318.9 | **failed** | 276 |
| `graphene` | 268 | 272 | **84*** | **84.3*** | 237 | 143 |

## 6.2 Performance comparison

In the remaining of this section, we use the convergence criterion (9) in all implementations of the LOBPCG algorithm to test the performance of these implementations. In addition to testing the conservative implementation of Hetmaniuk–Lehoucq algorithm (CHL), we also test an aggressive implementation (AHL) using the implicit updates discussed in Section 5.2 and assuming (12) (but not (13) and (14)). AHL performs fewer sparse matrix–vector multiplications (SpMV) at the expense of performing more implicitly updates. When SpMV is very cheap, which is the case here, this can lead to a slight increase in the overall time. In general, we expect SpMV to be more expensive in real world applications. Since the blockwise basis update strategy discussed in Section 5.1 is straightforward to include in both CHL and AHL, we turn this optimization on by default for these two solvers.

In Table 3, we report the number of iterations used by different implementations of the LOBPCG algorithm to reach convergence. For standard eigenvalue problems, all implementations use approximately the same number of iterations when they converge. One exception is `C60`, for which Anasazi takes fewer iterations than other solvers. Blopex terminates after 54 iterations. However 35 out of 176 computed eigenpairs are incorrect. Thus we report the number of iterations used by Blopex with $*$ for this problem. Both Anasazi and Blopex fail to converge within 2,000 iterations for `c-65` and `Ga10As10H30`. Anasazi also fails to converge for `Ga3As3H12`. For generalized eigenvalue problems, Anasazi takes more iterations to converge than our implementations of the LOBPCG algorithm. Blopex misses some of the desired eigenpairs for both the `nanotube` and `graphene` problems. Therefore, we mark the number of iterations for Blopex by $*$ for these problems in Table 3.

It is clear from this table that our implementations are more robust than Anasazi and Blopex for these examples, thanks to the basis selection strategies we use in our implementations.

Table 3: Number of iterations used by different implementations of the LOBPCG algorithm. Failures (including incorrect solutions marked with $*$) are highlighted in boldface.

| case | matrix | Anasazi | Blopex | AHL | CHL | Alg. 8 |
|---|---|---|---|---|---|---|
| 1 | `C60` | 39 | **54*** | 54 | 54 | 54 |
| 2 | `Si5H12` | 59 | 58 | 58 | 58 | 58 |
| 3 | `c-65` | **failed** | **failed** | 100 | 99 | 102 |
| 4 | `Andrews` | 52 | 56 | 55 | 55 | 55 |
| 5 | `Ga3As3H12` | **failed** | 63 | 65 | 65 | 65 |
| 6 | `Ga10As10H30` | **failed** | **failed** | 102 | 102 | 104 |
| 7 | nanotube | 208 | **446*** | 139 | 139 | 139 |
| 8 | graphene | 89 | **84*** | 52 | 52 | 52 |

We report the wall clock time used by different implementations of the LOBPCG algorithm on Edison in Figure 1. We plot the relative time, which is defined to be the ratio of the measured wall clock time used by a particular implementation of the LOBPCG algorithm over the wall clock time consumed by the CHL version of the LOBPCG implementation. We chose the CHL version as the baseline because it is a reliable implementation. We measured the wall clock time in both sequential runs and parallel runs with 24 cores. Because Anasazi and Blopex fail to converge within 2,000 iterations for some problems (e.g., `c-65`, `Ga3As3H12` and `Ga10As10H30`) or produce incorrect solutions for some problems (e.g. `C60`, `nanotube` and `graphene`), we report $2.5\times$ relative time for these tests in the figure.

With the exception of `Ga3As3H12`, for which Blopex runs slightly faster, the implementation of Algorithm 8 is the fastest among all implementation of the LOBPCG algorithm. It seems to be quite robust even though it skips some orthogonalization steps. It can be seen that both CHL and AHL take a bit more time than that used by Blopex. The additional time used in CHL and AHL is primarily due to the extra basis orthogonalization steps performed to improve the numerical stability of the LOBPCG algorithm.

In Figure 2, we use the `Andrews` problem as an example to show the breakdown of wall clock time spent in different components of several variants of the LOBPCG implementation. We can see that using the improved basis selection techniques and other performance optimization techniques in basis orthogonalization allows us to reduce the orthogonalization cost significantly in the implementation labeled by 4 in the figure.

When combined with the strategy of skipping some orthogonalization, the orthogonalization part of the computation is reduced from one third of the cost in the CHL and AHL versions to less than 10% of the total cost in version 5 of the implementation.

On average, CHL and AHL both require 3.4 iterations of `svqb`() in each call of `ortho`(), while Algorithm 8 reduces this number to 2.7. If we use one step of reorthogonalization (i.e., two iterations of `svqb`() per `ortho`() call) as the reference point, our new basis selection strategy reduces the additional reorthogonalization steps roughly by half.

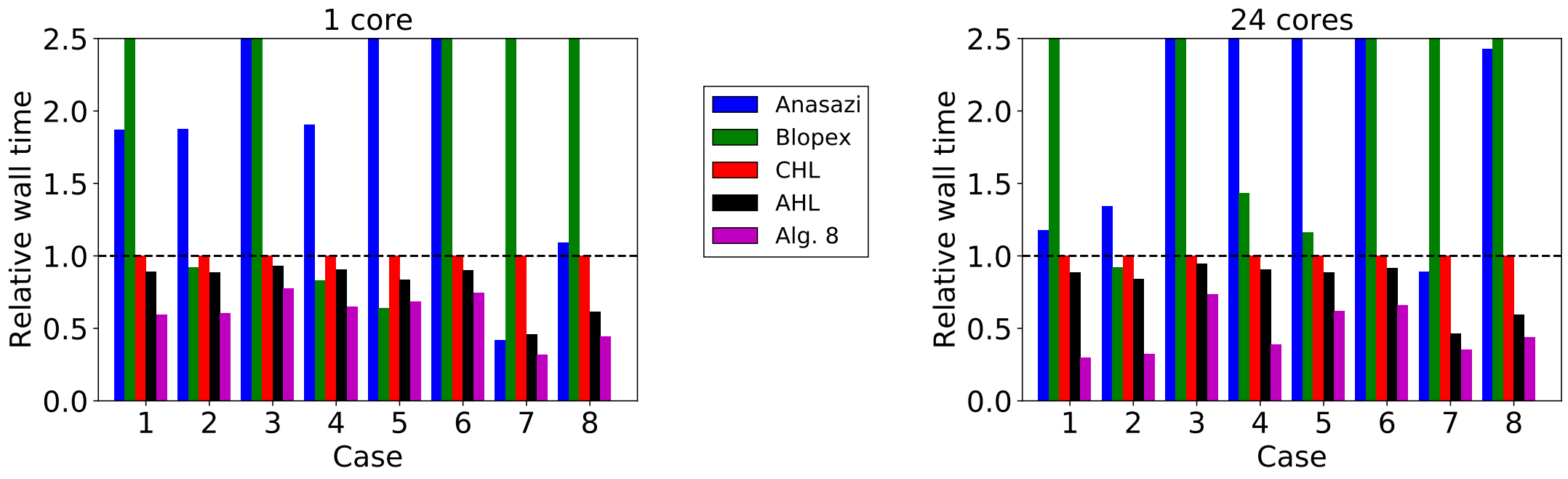

Figure 1: Execution time relative to the CHL implementation for tests performed with 1 core (left) and 24 cores (right). Failures (including incorrect solutions) are reported as 2.5× relative time.

The blockwise basis update and other techniques used in Algorithm 8 also make the orthogonalization cost nearly negligible.

The implicit update technique discussed in Section 5.2 is another consequence of the new basis selection strategy. If this technique is not adopted, the execution time for generating the projected matrix will increase by about 80% in Algorithm 8. Since this part of computation takes 10%–25% of the overall execution time, using our technique of implicit update leads to 5%–15% improvement in execution time.

Overall, the strategies discussed in Section 5 largely improve the performance of the LOBPCG algorithm, and the ones discussed in Section 4 do not cause additional overhead. A similar breakdown of wall clock time and the relative contributions of different techniques towards performance improvement of the LOBPCG implementation is observed for other test problems as well.

### 6.3 Scalability tests

Our implementation of Algorithm 8 demonstrates promising scalability on 24 cores in Figure 1. To have a closer look on the thread scalability, we perform tests on `Si5H12`, for which all implementations succeed, using different numbers of cores. The execution time is plotted in Figure 3(a).

The parallel scaling tests for Anasazi, AHL, and—to a lesser extent—CHL and Blopex, demonstrate an interesting artifact: Performance drops for twelve core tests. Note that each 24-core node on Edison is separated into two 12-core non-uniform memory access (NUMA) nodes. However the performance drop is not fully due to the NUMA penalty, as it does not occur when more than twelve cores are used. The implementation difference between CHL and AHL illuminates the source of this phenomenon. AHL performs additional tall skinny GEMM during implicit basis updates. This reveals an inefficiency in matrix multiply when the NUMA node saturates. The strategy of updating basis in a block row manner discussed in Section 5.1 requires a single synchronization—communicating the small update matrix—before becoming naively parallel. As a result, the phenomenon vanishes entirely in our implementation of Algorithm 8. Figure 3(b) shows the execution time spent on these tall skinny GEMM operations, which account for a significant portion (typically, 7%–25%) of the overall execution time. They have similar scalability compared to that the overall execution time. Hence, the plots confirm that these operations are indeed the source of this phenomenon.

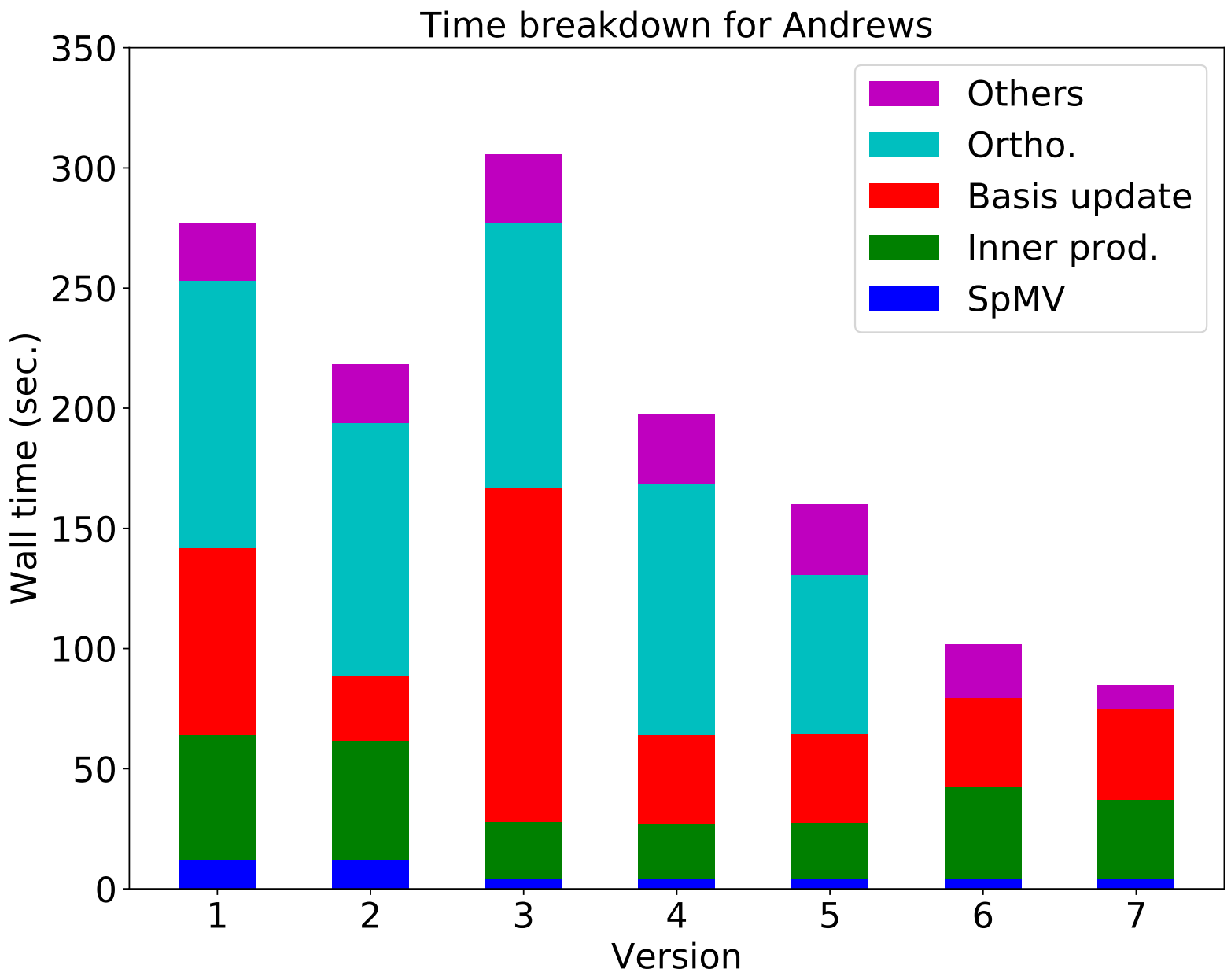

Figure 2: A breakdown of the wall clock time spent in different components of several versions of the LOBPCG implementation applied to the `Andrews` matrix. Each version is labeled by an integer with 1 being the CHL version *without* blockwise basis update, 2 being the CHL version (with blockwise basis update on top of version 1), 3 being the AHL version without blockwise basis update which performs implicit product updates on top of version 1, 4 being the AHL version (with blockwise basis update on top of version 3), 5 being a version that performs basis selection/basis trunction and other optimizations inside orthogonalization on top of version 4, 6 being a version that skips orthogonalization on top of version 5, and 7 being the version that implements Algorithm 8 that makes use of implicit product update (i.e., (13) and (14)) on top of version 6.

Finally we examine the effect of scaling the number of desired eigenpairs and the corresponding block dimension, because handling a large number of eigenpairs is one of the motivations of this work. We choose `Andrews` as the test case, as all of the implementations we test are successful on this matrix. The execution time on 24 cores is shown in Figure 4. Our implementation of Algorithm 8 consistently outperforms other implementations.

# 7 Concluding remarks

We developed a number of techniques to improve the stability of the LOBPCG algorithm, especially when the algorithm is used to compute a relatively large number of eigenpairs. We showed that a careful implementation of these techniques can also lead to improvement in the efficiency of the algorithm. We demonstrated the improvement by several numerical examples performed on single processor and multi-core systems using OpenMP parallelization. Most of the techniques discussed here can help improve the stability and performance of a distributed memory parallel implementation of the LOBPCG algorithm also (see, e.g., [14] for a recent implementation on distributed memory systems). However, the performance of such an implementation often depends on the parallelization of the sparse matrix–vector multiplications (SpMV). The best strategy for producing a scalable parallel SpMV is often problem dependent and is beyond the scope of this paper. We are in the process of developing such an implementation for general sparse symmetric matrices based on some recent techniques we developed for optimizing the performance of SpMV. We will demonstrate the performance of such implementation in a future publication.

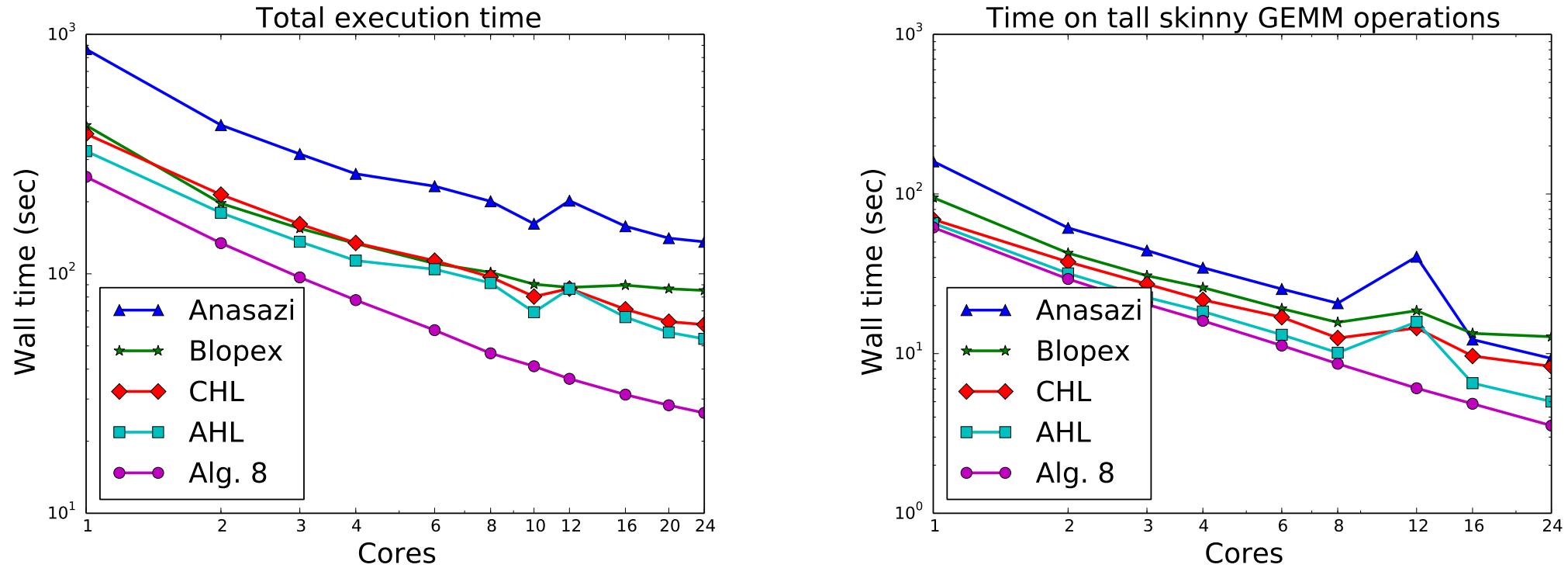

(a) (b)

Figure 3: Execution time for the test case `Si5H12` with respect to the number of cores: (a) total execution time; (b) time spent on tall skinny GEMM operations.

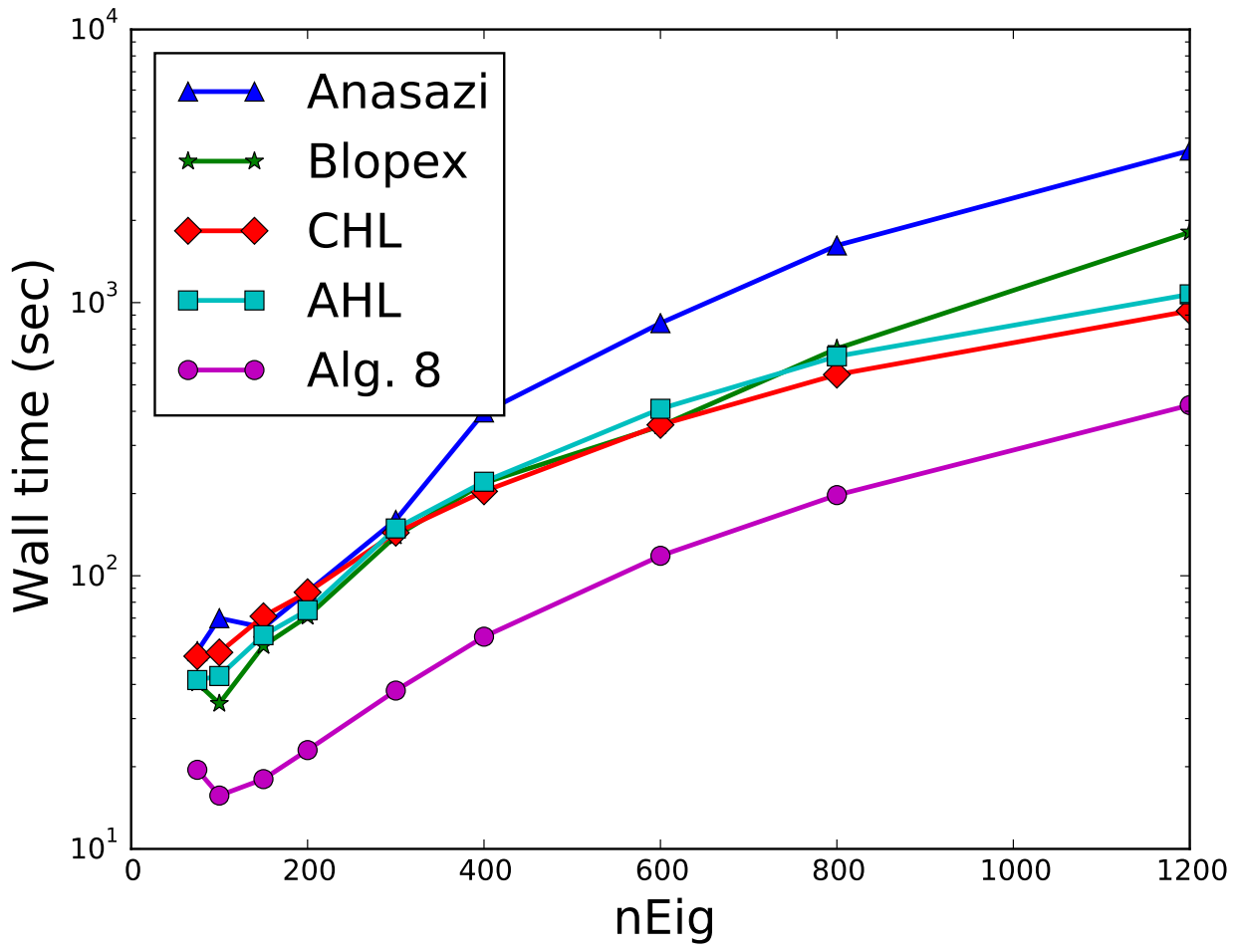

Figure 4: Execution time for the test case `Andrews` with respect to the number of desired eigenpairs.

## Acknowledgments

We would like to thank Lin Lin, Osni Marques, and Eugene Vecharynski for several useful discussions regarding this work. We thank Wei Hu for providing matrices produced from the SIESTA software for generalized eigenvalue problems.